\documentclass[a4paper,11pt]{article}
\usepackage{graphicx}
\usepackage{amsfonts}
\usepackage{amsmath}
\usepackage{amssymb}
\usepackage{indentfirst}
\usepackage{titlesec}
\usepackage{caption2}
\usepackage{color}

\def\barr{\begin{array}}
\def\earr{\end{array}}
\def\bali{\begin{aligned}}
\def\eali{\end{aligned}}
\def\bearr{\begin{eqnarray}}
\def\eearr{\end{eqnarray}}

\providecommand{\play}{\displaystyle}

\providecommand{\li}{\limits}

\providecommand{\pt}{\partial}

\providecommand{\ra}{\rightarrow}

\providecommand{\da}{\downarrow}

\providecommand{\Prob}{\mathbf P}

\providecommand{\E}{\mathbf E}

\providecommand{\al}{\alpha}

\providecommand{\bt}{\beta}

\providecommand{\gm}{\gamma}

\providecommand{\dt}{\delta}

\providecommand{\Dt}{\Delta}

\providecommand{\ve}{\varepsilon}

\providecommand{\tht}{\theta}

\providecommand{\kp}{\kappa}

\providecommand{\lb}{\lambda}

\providecommand{\sm}{\sigma}

\providecommand{\zt}{\zeta}

\providecommand{\R}{\mathbb R}

\providecommand{\cB}{\mathcal B}

\providecommand{\cE}{\mathcal E}

\providecommand{\cF}{\mathcal F}

\providecommand{\cL}{\mathcal L}

\providecommand{\cS}{\mathcal S}

\providecommand{\iY}{\mathfrak{Y}}

\providecommand{\iiY}{\mathbf{\mathfrak{Y}}}

\providecommand{\1}{\mathbf 1}

\providecommand{\contfunc}{\mathbf{C}}

\providecommand{\grad}{\nabla}

\providecommand{\vphi}{\varphi}

\providecommand{\graph}{\mathbb{G}}

\providecommand{\aone}{a^{(1)}}

\providecommand{\azero}{a^{(0)}}

\providecommand{\aonebar}{\overline{a^{(1)}}}

\providecommand{\Vol}{\text{Volume}}

\providecommand{\gmin}{\overline{\gm}}

\providecommand{\gmout}{\underline{\gm}}

\providecommand{\gmoout}{\underline{\underline{\gm}}}

\providecommand{\smtl}{\widetilde{\sm}}

\providecommand{\ntl}{\widetilde{n}}

\providecommand{\e}{\mathbf{e}}

\begin{document}

\title{On second order elliptic equations with a small parameter}
\author{Mark Freidlin\thanks{Department of Mathematics,
University of Maryland at College Park, mif@math.umd.edu.} \ , \
Wenqing Hu\thanks{Department of Mathematics, University of Maryland
at College Park, huwenqing@math.umd.edu.}}

\date{}

\maketitle

\begin{abstract}

The Neumann problem with a small parameter
$$\left(\dfrac{1}{\ve}L_0+L_1\right)u^\ve(x)=f(x) \text{ for } x\in G , \
\left.\dfrac{\pt u^\ve}{\pt \gm^\ve}(x)\right|_{\pt G}=0$$ is
considered in this paper. The operators $L_0$ and $L_1$ are
self-adjoint second order operators. We assume that $L_0$ has a
non-negative characteristic form and $L_1$ is strictly elliptic. The
reflection is with respect to inward co-normal unit vector
$\gm^\ve(x)$. The behavior of $\lim\li_{\ve\da 0}u^\ve(x)$ is
effectively described via the solution of an ordinary differential
equation on a tree. We calculate the differential operators inside
the edges of this tree and the gluing condition at the root. Our
approach is based on an analysis of the corresponding diffusion
processes.

\end{abstract}



\section{Introduction}

Let $G$ be a bounded domain in $\R^d$ ($d\geq 2$) with the smooth
boundary $\pt G$, $$L_k u(x)=\dfrac{1}{2}\sum\li_{i,j=1}^d
a^{(k)}_{ij}(x)\dfrac{\pt^2 u}{\pt x_i \pt x_j}+\sum\li_{i=1}^d
b_i^{(k)}(x)\dfrac{\pt u}{\pt x_i} \ , \ k=0,1 \ , \ x\in \R^d \ .$$

The coefficients are assumed to be smooth enough, say, in
$\contfunc^{(2)}(\R^d)$, i.e., having continuous second derivatives.

Boundary problems for the operator $L_\ve=L_0+\ve L_1$ in the domain
$G$ and corresponding initial-boundary problems for the equation
$\dfrac{\pt u^\ve(t,x)}{\pt t}=L_\ve u^\ve$, $t>0$, $x\in G$, are
considered. The operator $L_\ve$ is assumed to be elliptic for
$\ve>0$. One can study the limiting behavior of solutions of
stationary problems as $\ve\da 0$ and the limiting behavior of
solutions of initial-boundary problems as $\ve \da 0$ and $t\ra
\infty$.

If the operator $L_0$ is elliptic in $G\cup \pt G$, the problem is
simple: $u^\ve$ converges to the solution of corresponding problem
for the operator $L_0$. In the case of degenerate operator $L_0$,
situation is more complicated, and the question was considered in
numerous papers. First, the case of first order operator $L_0$ was
considered: $L_0=b^{(0)}(x)\cdot \grad$,
$b^{(0)}(x)=(b^{(0)}_1(x),...,b^{(0)}_d(x))$. N.Levinson
\cite{[Levinson]} showed in 1950-th that, if the characteristics of
$L_0$ (e.g., trajectories of the dynamical system
$\dot{X}_t=b^{(0)}(X_t)$ in $\R^d$) leave the domain $G$ in finite
time and cross the boundary in a regular way, then the solution of
the Dirichlet problem $L_\ve u^\ve=0$, $x\in G$, $u^\ve(x)|_{\pt
G}=\psi(x)$, converges as $\ve\da 0$ to the solution of degenerate
equation $L_0u^0(x)=0$, $x\in G$, with the boundary condition
$\psi(x)$ ($\psi(x)$ is assumed to be continuous) on the part of
$\pt G$ through which the characteristics leave the domain. Such a
solution $u^0(x)$ is unique.

Most of subsequent results concerning this problem were obtained by
probabilistic methods. With each operator $L_\ve$, $\ve\geq 0$, one
can (see \cite{[F Lemma]}, notice that the coefficients of
$a_{ij}^{(k)}(x)$ are in $\contfunc^{(2)}(\R^d)$) connect a
diffusion process $\widetilde{X}_t^\ve$ in $\R^d$ defined by the
stochastic differential equation

$$
\begin{array}{l}
\dot{\widetilde{X}}_t^\ve=b^{(0)}(\widetilde{X}_t^\ve)+\ve
b^{(1)}(\widetilde{X}_t^\ve)+
\sm^{(0)}(\widetilde{X}_t^\ve)\dot{W}_t^0+\sqrt{\ve}\sm^{(1)}(\widetilde{X}_t^\ve)\dot{W}_t^1
\ \ ,
\\
\widetilde{X}_0^\ve=x\in \R^d \ , \ t>0 \ , \
\sm^{(k)}(x)(\sm^{(k)}(x))^*=(a_{ij}^{(k)}(x))=a^{(k)}(x) \ , \ k=0
, 1 \ .
\end{array}$$

Here $W_t^0$ and $W_t^1$ are independent Wiener processes in $\R^d$.
Then the solution of the Dirichlet problem for the equation $L_\ve
u^\ve(x)=0$, $x\in G$, and of the initial boundary problem for
$\dfrac{\pt u^\ve(t,x)}{\pt t}=L_\ve u^\ve(t,x)$ can be represented
as expectations of corresponding functionals of
$\widetilde{X}_t^\ve$. The trajectories $\widetilde{X}_t^\ve$, in a
sense, play the same role as characteristics in the case of first
order operator $L_0$. Using these representations and studying
limiting behavior of process $\widetilde{X}_t^\ve$ one can describe
the limiting behavior of the boundary problems (see \cite{[F red
book]}, \cite{[FW book]}).

If problems with the Neumann boundary conditions are considered, one
can use the corresponding diffusion process with reflection on the
boundary (see, for instance, \cite[\S 2.5]{[F red book]}). Various
cases of first order operators $L_0$ not satisfying Levinson's
conditions were examined using the probabilistic approach (see
\cite{[F red book]}, \cite{[FW book]} and the references therein).

If the operator $L_0$ has terms with second derivatives, one can
introduce a generalized Levinson condition (\cite[\S 4.2]{[F red
book]}). Under this condition the equation $L_0 u^0(x)=0$, $x\in G$,
with appropriate Dirichlet type boundary conditions has a unique
solution, and the solution $u^\ve(x)$ of the Dirichlet problem for
equation $L_\ve u^\ve(x)=0$, $x\in G$ converges to $u^0(x)$ as $\ve
\da 0$. The difference with the classical Levinson case is just in
the rate of convergence: under mild additional assumptions
$|u^\ve(x)-u^0(x)|<\ve^\gm$ for some $\gm>0$ and $0<\ve<\!<1$, but
for any $\gm'>0$ one can find $L_0$ with infinitely differentiable
coefficients non-degenerating on $\pt G$ such that
$|u^\ve(x)-u^0(x)|$ is greater than $\ve^{\gm'}$ at a point $x\in G$
and $0<\ve<\!<1$.

A convenient way to specify the degeneration of $L_0$ is given by
the conservation laws. A function $H(x)$ is called a first integral
for the process $X_t^0$ corresponding to $L_0$ if $\Prob_x(X_t^0\in
S(H(x)))=1$ for all $t\geq 0$ and $x\in \R^d$, where $S(z)=\{y\in
\R^d: H(y)=z\}$; here and below the subscript $x\in \R^d$ in the
probability $\Prob_x$ or expected value $\E_x$ means that the
trajectory of the process starts at the point $x$.

We consider in this paper self-adjoint operators $L_0$ and $L_1$:
$$L_k u(x)=\dfrac{1}{2}\grad \cdot (a^{(k)}(x)\grad u(x)) \ .$$ Then
a smooth function $H(x)$ is a first integral for the process
$\widetilde{X}_t^0$ (for the corresponding operator $L_0$) if and
only if $\azero(x)\grad H(x) \equiv 0$. In general, the process
$\widetilde{X}_t^0$ can have several independent smooth first
integrals. To restrict ourselves to the case of one smooth first
integral we assume that $\e\cdot(\azero(x)\e)\geq
\underline{a}(x)|\e|_{\R^d}^2$ for each $\e\in \R^d$ such that
$\e\cdot \grad H(x)=0$: It is assumed that $\underline{a}(x)$ is
smooth and strictly positive if $x$ is not a critical point of
$H(x)$; if $x_0$ is a critical point, $\azero(x_0)=0$ and
$\underline{a}(x_0)=0$.

To be specific we consider the Neumann problem

$$\dfrac{1}{\ve}L_\ve u^\ve=\left(\dfrac{1}{\ve}L_0+L_1\right)u^\ve(x)=f(x) \ , \
\left.\dfrac{\pt u^\ve}{\pt \gm^\ve}(x)\right|_{\pt G}=0 \ ;
\eqno(1.1)$$ $\gm^\ve(x)$ here is the inward co-normal unit vector
to $\pt G$ corresponding to $L_\ve$. Let $X_t^\ve$ be the process in
$G\cup \pt G$ governed by the operator inside $G$ with reflection
along the co-normal to $\pt G$. Since $L_\ve$ is self-adjoint, the
Lebesgue measure is invariant for the process $X_t^\ve$, and the
problem (1.1) is solvable for continuous $f(x)$ such that
$\play{\int_{G}f(x)dx=0}$. Together with the last condition, we
assume that $L_1$ is not degenerate in $G \cup \pt G$, so that to
single out a unique solution of (1.1), we shall fix the value of
$u^\ve(x)$ at a point $x_{O}\in G\cup \pt G$ which is fixed the same
for all $\ve>0$. We let $u^\ve(x_O)=0$.

Then the solution of problem (1.1) can be written in the form (see,
for instance \cite{[F red book]})

$$u^\ve(x)=-\int_0^\infty \E_x f(X_t^\ve)dt+\int_0^\infty \E_{x_O} f(X_t^\ve)dt  \ . \eqno(1.2)$$

If the first integral $H(x)$ has in $G \cup \pt G$ no critical
points, one can describe the $\lim\li_{\ve\da 0}u^\ve(x)$ in the way
similar to \cite{[FW fish]}: One shall introduce a graph $\graph$
corresponding to the set of connected components of the
intersections of the level sets of $H(x)$ within $G$. A boundary
problem on $\graph$ with appropriate gluing conditions at the
vertices can be formulated, and the solution of this problem defines
$\lim\li_{\ve\da 0}u^\ve(x)$. If the function $H(x)$ has saddle
points inside $G$, additional branchings in the graph appear. The
gluing conditions at these new vertices can be calculated using the
results of \cite{[FWeb]}.

All mentioned above results concern the case when the rank of
$\azero(x)$ is constant and equal to $d-1$ for all $x\in G\cup \pt
G$ except the critical points of $H(x)$. In this paper, we consider
the case when $L_0$ is non-degenerate in a connected subdomain
$\cE\subset G$, and we let $H(x)$ be equal to a constant on $\cE$.
Outside $\cE$ the first integral $H(x)$ has a finite number of
critical points (see Fig.1). For convenience of presentation, we
shall then introduce several first integrals $H_k(x)$ ($k=1,...,r$)
for each of the connected components $U_1,...,U_r$ on which $L_0$ is
degenerate. We shall let $H(x)=H_k(x)$ for $x\in U_k$. A more
concrete setup of the problem is in Section 2. Existence of the
domain $\cE$ where the operator $L_0$ is not degenerate leads to
more general gluing conditions. The limiting process on the graph
spends a positive time at the vertex corresponding to $\cE$.

\begin{figure}
\renewcommand{\captionlabeldelim}{.}
\centering
\includegraphics[height=5cm, width=14cm , bb=11 82 399 217]{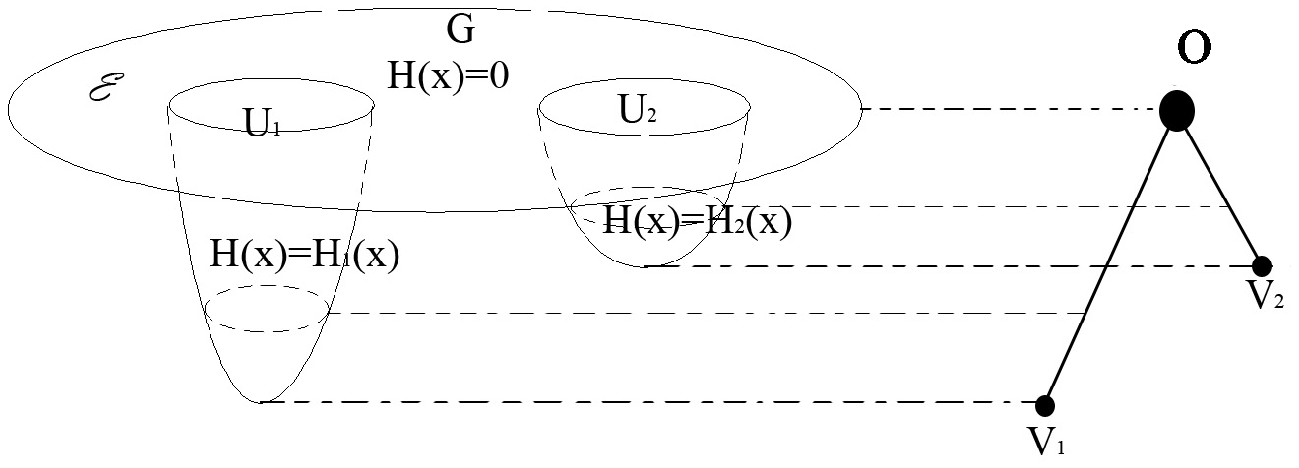}
\caption{}
\end{figure}

Let $S(z)=\{x\in G\cup \pt G: H(x)=z\}$. The graph $\graph$ is the
result of identification of points of each connected component of
every level set $S(z)$. Let $\iiY: G\cup \pt G \ra \graph$ be the
identification mapping. We call $\iiY(x)$ the projection of $x$ onto
$\graph$. We consider the projection $Y_t^\ve=\iiY(X_t^\ve)$ of the
process $X_t^\ve$ on $\graph$ and prove that processes $Y_t^\ve$ on
$\graph$ converge weakly in the space of continuous functions $[0,T]
\ra \graph$ to a diffusion process $Y_t$ on $\graph$. The process
$Y_t$ is defined by a family of differential operators, one on each
edge, and by gluing conditions at the vertices. We calculate the
operators and the gluing condition. The function
$u^0(x)=\lim\li_{\ve\da 0}u^\ve(x)$ is constant on each connected
component of every level set of $H(x)$: $u^0(x)=v(\iiY(x))$. We
formulate a boundary problem for the function $v(y)$, $y\in \graph$,
which has a unique solution, and actually can be solved explicitly.

The organization of this paper is as follows: Section 2 sets up the
problem and gives the main results. Section 3 is devoted to the
proof of the main results in Section 2. Section 4 proves auxiliary
results needed in Section 3.

\section{Main results}

Let us first speak about our assumptions.

Suppose we have a bounded domain $G\subset \R^d$, with smooth
boundary $\pt G$. We assume that $d\geq 2$. Let $L_0$ be a
self-adjoint operator
$$L_0 u (x)=\dfrac{1}{2}\sum\li_{i,j=1}^d \dfrac{\pt}{\pt x_i}\left(\azero_{ij}(x) \dfrac{\pt u(x)}{\pt
x_j}\right)=\dfrac{1}{2}\grad\cdot (\azero(x)\grad u(x)) \ .$$

Let $U_1,...,U_r$ be several regions inside $G$. They are simply
connected open sets and their closure do not intersect each other.
Let us assume that the matrix
$\azero(x)\equiv(\azero_{ij}(x))_{1\leq i,j \leq d}$ is positive
definite on $\cE=[G]\backslash \left(\cup_{k=1}^r [U_k]\right)$
(here $[D]$ is the closure of a domain $D$). For $x\in [\cE]$, the
coefficients $\azero_{ij}(x), 1\leq i,j \leq d$ are assumed to be in
$\contfunc^{(3)}([\cE])$.

Let us discuss the case when $x\in \cup_{k=1}^r[U_k]$. For each
$k=1,...,r$ and $x\in [U_k]$, the coefficients $\azero_{ij}(x),
1\leq i,j \leq d$ are assumed to be in $\contfunc^{(3)}([U_k])$. We
assume that the matrix $(\azero_{ij}(x))_{1\leq i,j \leq d}$ is
degenerate on $\cup_{k=1}^r[U_k]$. To specify this degeneration, we
assume that within each $[U_k]$ there is a function $H_k$ which is a
first integral of the (degenerate) operator $L_0$, i.e.,
$\azero(x)\grad H_k(x)=0$ for $x\in [U_k]$. The function $H_k$ is
assumed to be in at least $\contfunc^{(4)}([U_k])$. Let $H_k$ have
only one minimum $m_k$ inside $U_k$. (We can always make this
assumption since if $m_k$ is a maximum we work with $-H_k$ instead
of $H_k$.) Let $x_k(m_k)$ be the point in $U_k$ corresponding to the
minimum $m_k$. This minimum is assumed to be non-degenerate, i.e.,
the matrix $\left(\dfrac{\pt^2 H_k}{\pt x_i\pt
x_j}(x_k(m_k))\right)_{1\leq i,j \leq d}$ is positive definite.
Since the choice of $H_k$ is up to a constant we can assume that
$H_k=0$ on $\pt U_k$. For $h\in (m_k,0]$ the level surfaces
$C_k(h)=\{x\in U_k: H_k(x)=h\}$ of the functions $H_k$ inside $U_k$
are closed surfaces of dimension $(d-1)$ and the operator $L_0$ is
non-degenerate on $C_k(h)$. Let $\gm_k=\pt U_k=C_k(0)$. A
non-degeneracy condition of $\azero(x)$ on $C_k(h)$ is assumed: for
any vector $\mathbf{e}\in \R^d$ such that $\mathbf{e}\cdot \grad
H_k=0$ we have $\mathbf{e}\cdot(\azero(x) \mathbf{e})\geq
\underline{a}(x)|\mathbf{e}|_{\R^d}^2$ for some $\underline{a}(x)>0$
and $x\neq x_k(m_k)$. We set $\azero(x_k(m_k))=0$ and
$C_k(m_k)=\{x_k(m_k)\}$. We assume that the level surfaces $C_k(h)$
for $h\in (m_k,0]$ divide $U_k\backslash \{x_k(m_k)\}$ into pieces
of $(d-1)$-dimensional surfaces.

For simplicity of presentation we will assume that $\grad H_k
(x)\neq 0$ for $x\in \gm_k$. One can introduce a global first
integral $H(x)$ on $[G]$ as in Section 1: $H(x)=H_k(x)$ for $x\in
U_k$ and $H(x)=0$ for $x\in [\cE]$. We notice that the function
$H(x)$ so defined as a global first integral is not smooth at
$\cup_{k=1}^r \gm_k$. However, this will not affect our analysis.

Let $\gm=\cup_{k=1}^r \gm_k$. We will define a unit vector field
$\e_d(x)$ in a small neighborhood of $\gm_k$ at the beginning of
Section 4. Roughly speaking, the vector field $\e_d(x)$ is the
direction in which the matrix $\azero(x)$ degenerates. We assume
that the order of degeneracy is given by the condition that for this
$\e_d(x)$ as $x$ belongs to the intersection of a small neighborhood
of $\cup_{k=1}^r \gm_k$ and $\cE$ we have
$$\text{const}_1\cdot\text{dist}^2(x,\gm)\leq
\e_d(x)\cdot(\azero(x)\e_d(x))\leq
\text{const}_2\cdot\text{dist}^2(x,\gm)$$ for some $\text{const}_1$,
$\text{const}_2>0$. The distance $\text{dist}(x,\gm)$ is the
Euclidean distance between $x$ and $\gm$. The vector field
$\e_d(x)=\dfrac{\grad H_k}{|\grad H_k|_{\R^d}}$ for $x\in \gm_k$.

In particular, our assumptions imply that the matrix $\azero(x)$ has
rank $d$ in $\cE$ and rank $(d-1)$ in $\cup_{k=1}^r ([U_k]\backslash
\{x_k(m_k)\})$. At the points $x_k(m_k)$ the matrix $\azero(x)$ is
just $0$. However, the coefficients $\azero_{ij}(x)$, $1\leq i,j\leq
d$ are only in $\contfunc^{(1)}$ for $x\in \gm$. We notice that in
this case results of \cite{[F Lemma]} do not apply. We shall then
assume that there is a decomposition
$\azero(x)=\sm^{(0)}(x)(\sm^{(0)}(x))^*$ for all $x\in [G]$. The
square matrix $\sm^{(0)}(x)$ has bounded Lipschitz continuous terms.

We shall assume that the operator $L_1$ governing the perturbation
is self-adjoint and strictly elliptic within $[G]$: $$L_1 u
(x)=\dfrac{1}{2}\sum\li_{i,j=1}^{d}\dfrac{\pt}{\pt
x_i}\left(\aone_{ij}(x)\dfrac{\pt u (x)}{\pt
x_j}\right)=\dfrac{1}{2}\grad\cdot (\aone(x)\grad u(x)) \ .$$

Again we denote the matrix $\aone(x)\equiv(\aone_{ij}(x))_{1\leq i,j
\leq d}$ and we assume that the terms $\aone_{ij}(x)$ are in class
$\contfunc^{(2)}(\R^d)$. In this case results of \cite{[F Lemma]}
apply and we have $\aone(x)=\sm^{(1)}(x)(\sm^{(1)}(x))^*$ for all
$x\in [G]$. The square matrix $\sm^{(1)}(x)$ have bounded Lipschitz
continuous terms.

Let us put a Neumann boundary condition with respect to co-normal
unit vector $\gm^\ve(x)$ pointing inward on $\pt G$. Let $X_t^\ve$
be the diffusion process in $[G]$, corresponding to the operator
$\dfrac{1}{\ve}L_0+L_1$ inside $G$ with co-normal reflection at $\pt
G$. We see that Lebesgue measure is invariant for the process
$X_t^\ve$.

\

Consider a simple example. We use polar coordinates for the plane
$\R^2$. Let $G=\{(r,\tht): 0\leq \tht\leq 2\pi, 0\leq r < 2\}$. Let
$U_1=\{(r,\tht): 0\leq \tht \leq 2\pi, 0\leq r <1\}$. Let
$H_1(r,\tht)=\dfrac{r^2-1}{2}$. Let $H(r,\tht)=H_1(r,\tht)$ for
$0\leq r<1$ and $H(r,\tht)=0$ for $1\leq r\leq 2$. Let
$\lb(r)=(r-1)^2$ for $1\leq r\leq 2$ and $\lb(r)=0$ for $0\leq r<
1$. Let
$$a^{(0)}(r,\tht)=\begin{pmatrix}\lb(r)\cos^2\tht+r\sin^2\tht & (\lb(r)-r)\sin\tht\cos\tht\\
(\lb(r)-r)\sin\tht\cos\tht &
\lb(r)\sin^2\tht+r\cos^2\tht\end{pmatrix} \ .$$ Let $\e_r=(\cos\tht,
\sin\tht)$ and $\e_\tht=(-\sin\tht, \cos\tht)$. Then it is easy to
check that $a^{(0)}(r,\tht)\e_r=\lb(r)\e_r$ and
$a^{(0)}(r,\tht)\e_\tht=r\e_\tht$. We have the decomposition
$a^{(0)}(x)=\sm^{(0)}(x)(\sm^{(0)}(x))^*$ where
$$\sm^{(0)}(r,\tht)=\begin{pmatrix}\sqrt{\lb(r)}\cos^2\tht+\sqrt{r}\sin^2\tht & (\sqrt{\lb(r)}-\sqrt{r})\sin\tht\cos\tht\\
(\sqrt{\lb(r)}-\sqrt{r})\sin\tht\cos\tht &
\sqrt{\lb(r)}\sin^2\tht+\sqrt{r}\cos^2\tht\end{pmatrix} \ .$$

\

Let us then speak about the results.

We construct a graph $\graph$ as follows. The graph $\graph$ has $r$
edges $I_1,...,I_r$ joined together at one vertex $O$. Let the other
endpoint of $I_k$ be $V_k$. Let us write $I_k=[m_k,0]$. The
coordinate $(k,H_k)$ is a global coordinate on $\graph$. The root
$O$ corresponds to all $(k,0)$ for $k=1,...,r$. Let us introduce an
identification map $\iY: [G] \ra \graph$: for $x\in [\cE]$ we have
$\iY(x)=O$ and for $x\in U_k$ we have $\iY(x)=(k,H_k(x))$. Let the
process $Y_t^\ve=\iY(X_t^\ve)$. We are going to prove, that as $\ve
\da 0$ the processes $Y_t^\ve$ converge weakly in the space
$\contfunc_{[0,T]}(\graph)$ to a Markov process $Y_t$ on $\graph$.

The process $Y_t$ is defined as follows. It is a diffusion process
on the graph $\graph$ with generator $A$ and the domain of
definition $D(A)$. Inside each $I_k$ it is governed by an operator
$\cL_k$ defined as $$\cL_k f(k,H_k)=
\dfrac{1}{2}M_k^{-1}(H_k)\dfrac{d}{dH_k}\left(M_k(H_k)\aonebar
(H_k)\dfrac{df}{dH_k}\right)\ .$$ Here
$$\aonebar(h)=M_k^{-1}(h)\int\li_{C_k(h)}\dfrac{(\aone(x)\grad H_k(x),\grad H_k(x))}
{|\grad H_k(x)|_{\R^d}}d\sm \ ,$$ and normalizing factor
$$M_k(h)=\int\li_{C_k(h)}\dfrac{d\sm}{|\grad H_k(x)|_{\R^d}} \
.$$ The notation $d\sm$ denotes the integral with respect to the
area element on $C_k(h)$.

We set $Af=\cL_k f$ for $(k,H_k)\in (I_k)$ ($(I_k)$ is the interior
of the interval $I_k$). Let the limit $\lim\li_{(k,H_k)\ra
O}Af(k,H_k)$ be finite and independent of $k$. This limit is set to
be $Af(O)$. The domain of definition $D(A)$ of the operator $A$
consists of those functions $f$ that are twice continuously
differentiable inside each $I_k$ having the limit $\lim\li_{H_k\ra
0}\dfrac{\pt f}{\pt H_k}(k,H_k)$. These functions satisfy the gluing
condition at the vertex $O$:
$$0=\Vol(\cE)\cdot Af(O)+\dfrac{1}{2}\sum\li_{k=1}^r p_k \cdot \lim\li_{H_k\ra
0}\dfrac{\pt f}{\pt H_k}(k,H_k) \ . \eqno(2.1)$$ Here $\Vol(\cE)$ is
the $d$-dimensional volume of the domain $\cE$ and
$$p_k=\int_{\gm_k}\dfrac{(a^{(1)}(x)\grad H_k(x), \grad H_k(x))}{|\grad H_k(x)|_{\R^d}}d\sm \
.$$

For the exterior vertices $V_1,...,V_r$ no additional assumptions
are to be imposed on the behavior of the function $f$ in the domain
$D(A)$.

It was proved in \cite{[FW Diffusion process on a graph]} the the
process $Y_t$ exists and is a strong Markov process on the graph
$\graph$.

We have

\

\textbf{Theorem 2.1.} \textit{As $\ve \da 0$ the processes $Y_t^\ve$
converge weakly to $Y_t$ in $\contfunc_{[0,T]}(\graph)$.}

\

Let $\mu_x^\ve$ be the distribution of the trajectory
$Y_t^\ve=\iiY(X_t^\ve)$ starting from a point $x\in [G]$ in the
space $\contfunc_{[0,T]}(\graph)$: for each Borel subset $B\subseteq
\contfunc_{[0,T]}(\graph)$ we set
$\mu_x^\ve(B)=\Prob_{X_0^\ve=x}(Y_{\bullet}^\ve\in B)$. Similarly,
for each $y\in \graph$ we let $\mu_y^0$ be the distribution of $Y_t$
in the space $\contfunc_{[0,T]}(\graph)$ with
$\mu_y^0(B)=\Prob_{Y_0=y}(Y_{\bullet}\in B)$. Theorem 2.1 can be
reformulated as

\

\textbf{Theorem 2.2.} \textit{For every $x\in [G]$ and every $T>0$
the distribution $\mu_x^\ve$ converges weakly to $\mu_{\iiY(x)}^0$
as $\ve \da 0$. For every bounded continuous functional $F$ on
$\contfunc_{[0,T]}(\graph)$ we have $$\E_{X_0^\ve=x}
F(Y_{\bullet}^\ve)\ra \E_{Y_0=\iiY(x)}F(Y_{\bullet})$$ as $\ve \da
0$.}

\

The process $Y_t^\ve$ can be viewed as the slow component of the
process $X_t^\ve$. The fast component $Z_t^\ve$ of $X_t^\ve$ is a
process governed by the operator $\dfrac{1}{\ve}L_0$. The process
$Z_t^\ve$ moves on $\iiY^{-1}(y)$ for each $y\in \graph$: it is
moving on $[\cE]$ when $y=O$ and it is moving on $C_k(H_k)$ when
$y=(k,H_k)$. Since Lebesgue measure is invariant for the process
$X_t^\ve$, the fast component $Z_t^\ve$, as $\ve>0$ is small, has,
approximately, a distribution with density
$\dfrac{1}{\text{Volume}(\cE)}$ on $[\cE]$ (with respect to Lebesgue
measure on $\R^d$) and $\dfrac{1}{M_k(H_k)}\dfrac{1}{|\grad
H_k|_{\R^d}}$ on $C_k(H_k)$ (with respect to the area element $d\sm$
on $C_k(H_k))$. Using this we can formulate the above two theorems
in terms of differential equations:

\

\textbf{Theorem 2.3.} \textit{Consider the Neumann problem}

$$\dfrac{1}{\ve}L_\ve u^\ve(x)=\left(\dfrac{1}{\ve}L_0+L_1\right)u^\ve(x)=f(x) \textit{\text{ for }} x\in G \ ,
\ \left.\dfrac{\pt u^\ve}{\pt \gm^\ve}(x)\right|_{x\in \pt G}=0$$
\textit{with a H\"{o}lder continuous function $f(x)$ satisfying
$\play{\int_G f(x)dx=0}$. Let $u^\ve(x_O)=0$ for some $x_O\in G$.
Then we have }$$\lim\li_{\ve\da 0}u^\ve(x)=v(\iiY(x))$$\textit{
where $v(y)$ is a continuous function on $\graph$ such that }$$\cL_k
v(y)=-\bar{f}(y) \textit{\text{ for }} y\in (I_k) \ , \ k=1,...,r \
.$$\textit{ Here}
$$\bar{f}(y)=\dfrac{1}{\text{Volume}(\cE)}\int_{\cE}f(x)dx$$\textit{ when $y=O$
and}
$$\bar{f}(y)=\dfrac{1}{M_k(H_k)}\int_{C_k(H_k)}f(x)\dfrac{d\sm}{|\grad
H_k(x)|_{\R^d}}$$\textit{ when $y=(k,H_k)$.} \textit{ The function
$v(y)$ satisfies the gluing condition (2.1) and $v(\iiY(x_O))=0$.
Such a function $v(y)$ is unique. }

\section{Proof of Theorem 2.1}

The \textbf{Proof} of Theorem 2.1 follows the arguments of \cite{[FW
book]}, \cite{[FW AMS]}, \cite{[DK1]} and \cite{[DK2]}.

Heuristically, the idea of \cite{[DK2]} can be explained as follows.
The process $X_t^\ve$ moves within $[G]$ and has Lebesgue measure as
its invariant measure. Since the process $X_t^\ve$ has a "fast"
component governed by the operator $\dfrac{1}{\ve}L_0$, it will
spend a positive amount of time proportional to $\Vol(\cE)$ within
$\cE$ as $\ve\da 0$. As we project $X_t^\ve$ onto the graph $\graph$
and the whole ergodic component $\cE$ corresponds to $O$, the
limiting process $Y_t$ has a boundary condition with a "delay" at
$O$. (We recommend a nice article \cite{[Ito-McKean 1963]} and a
brief summary \cite[\S 5.7]{[Ito-McKean book]} about this boundary
condition.) Our gluing condition (2.1) ensures that the process
$Y_t$ has an invariant measure on $\graph$ that agrees with the
Lebesgure measure on $[G]$. We also refer to \cite[Ch.8, pp.
347--350]{[FW book]} for an explanation of this.

Let us first introduce some notations. Below we will often suppress
the small parameter $\ve$ and it could be understood directly from
the context. Let $\gmin_k=C_k(-\ve^{1/2})$ and $\gmin=\cup_{k=1}^r
\gmin_k$. Let $\sm$ be the first time when the process $X_t^\ve$
hits $\gm$. Let $\tau$ be the first time when the process $X_t^\ve$
hits $\gmin$. Let $\sm_0=\sm$. Let $\tau_n$ be the first time
following $\sm_n$ when the process reaches $\gmin$. For $n\geq 1$
let $\sm_n$ be the first time after $\tau_{n-1}$ when the process
$X_t^\ve$ reaches $\gm$.

Let $\sm^*\in \{\sm_0,\sm_1,...\}$ and we denote by $m_{\sm^*}^x$
the measure on $\gm$ induced by $X_{\sm^*}^\ve$ starting at $x$.
That is,

$$m_{\sm^*}^x(A)=\Prob_x(X_{\sm^*}^\ve\in A) \ , \ A\in \cB(\gm) \ .$$

Let $\nu(\bullet)$ be the invariant measure of the induced chain
$X_{\sm_n}^\ve$ on $\gm$. The key lemma of \cite{[DK2]} is the
following

\

\textbf{Lemma 3.1.} \textit{Let $x\in [\cE]$. For each $\dt>0$ and
all sufficiently small $\ve$ there is a stopping time $\sm^*$ which
may depend on $\dt, \ve$ and $x$ and such that $$\E_x\sm^*\leq \dt \
, \eqno(3.1)$$
$$\sup\li_{x\in \gm}\text{Var} (m_{\sm^*}^x(dy)-\nu(dy))\leq \dt \ , \eqno(3.2)$$
where $\text{Var}$ is the total variation of the signed measure. }

\

Our proof of this lemma is similar while a bit simpler than that of
\cite{[DK2]}.

\

\textbf{Proof.} We will prove, in Lemma 4.11 that $X_{\sm_n}^\ve$
satisfies the Doeblin condition on $\gm$ uniformly in $\ve$. This
implies that one can choose an $N$ depending on $\dt$ but
independent of $\ve$ such that the distribution of $X_{\sm_N}^\ve$
is $\dt$-close to the invariant measure $\nu(\bullet)$ on $\gm$.
That is, as we set $\sm^*=\sm_N$ the condition (3.2) is satisfied.

We are going to prove in Lemmas 4.8, 4.9, 4.10, respectively, that
$$\lim\li_{\ve\da 0}\sup\li_{x\in [\cE]} \E_x\sm=0  \ , \eqno(3.3)$$
$$\lim\li_{\ve\da 0}\sup\li_{x\in \gm} \E_x\tau=0  \ , \eqno(3.4)$$
$$\lim\li_{\ve\da 0}\sup\li_{x\in\gmin} \E_x\sm=0  \ . \eqno(3.5)$$

We can write
$\sm_N=\sum\li_{i=1}^N[(\sm_i-\tau_{i-1})+(\tau_{i-1}-\sm_{i-1})]+\sm_0$.
For each $i=1,...,N$ the random variable $\sm_i-\tau_{i-1}$ has the
same distribution as $\sm$ for the process $X_t^\ve$ starting at
some point on $\gmin$; similarly, the random variable
$\tau_{i-1}-\sm_{i-1}$ has the same distribution as $\tau$ for the
process $X_t^\ve$ starting at some point on $\gm$. The results
(3.3), (3.4), (3.5) imply that as $\ve$ is small (notice that $N$ is
fixed at this stage) we can choose $\sm^*=\sm_N$ and the condition
(3.1) is also satisfied. $\square$

\

\textbf{Proof of Theorem 2.1.} The proof of Theorem 2.1 is the same
as the proof of Lemma 2.1 (including the proof of Lemma 3.4) stated
in \cite{[DK2]} using the above Lemma 3.1. For the sake of
completeness let us briefly repeat it here. Reasoning as in
\cite{[DK1]}, \cite{[DK2]}, \cite{[FW book]}, \cite{[FW AMS]}, it
suffices to prove that for a function $f\in D(A)$, for every $T>0$
and uniformly in $x\in [G]$ we have

$$\E_x\left[f(H(X_T^\ve))-f(H(X_0^\ve))-\int_0^T Af(H(X_s^\ve))ds\right]\ra 0$$
as $\ve \da 0$. Here $H(x)=H_k(x)$ if $x\in U_k$ and $H(x)=0$ if
$x\in [\cE]$. Let us replace the time interval $[0,T]$ by a larger
one $[0,\smtl]$, where $\smtl$ is the first of the stopping times
$\sm_n$ which is greater than or equal to $T$: $\smtl=\min\li_{n:
\sm_n>T}\sm_n$. Let $\smtl=\sm_{\widetilde{n}+1}$. We have

$$\begin{array}{l}
\play{\E_x\left[f(H(X_T^\ve))-f(H(X_0^\ve))-\int_0^T
Af(H(X_s^\ve))ds\right]}
\\
=\play{\E_x\left[f(H(X_\sm^\ve))-f(H(X_0^\ve))-\int_0^{\sm}
Af(H(X_s^\ve))ds\right]+\E_x\sum\li_{k=0}^{\ntl}\int_{\sm_k}^{\sm_{k+1}}Af(H(X_s^\ve))ds-}
\\
\play{ \ \ \ \ \ \
-\E_x\E_{X_T^\ve}\left[f(H(X_{\sm}^\ve))-f(H(X_0^\ve))-\int_0^{\sm}Af(H(X_s^\ve))ds\right]}
\\
=(I)+(II)-(III) \ .
\end{array}$$

If $x\in \cup_{k=1}^r U_k$ we have $|(I)|\ra 0$ uniformly in $x$ as
$\ve \da 0$ due to averaging principle. If $x\in [\cE]$ then
$|(I)|\ra 0$ uniformly in $x$ as $\ve \da 0$ due to Lemma 4.8. In a
similar way we see that $|(III)|\ra 0$ uniformly in $x\in [G]$ as
$\ve \da 0$.

Let $\play{\al_k=\int_{\sm_k}^{\sm_{k+1}}Af(H(X_s^\ve))ds}$ and let
$\play{\bt_k=\sum\li_{n=0}^\infty \E_x(\al_{k+n}|\cF_k)}$ ($\cF_k$
is the filtration generated by the process $X_t^\ve$ for $t\leq
\sm_k$). We have $\E_x(\al_k-\bt_k+\bt_{k+1}|\cF_k)=0$ and therefore
$\play{\left(\sum\li_{k=1}^n
(\al_k-\bt_k+\bt_{k+1}),\cF_{n+1}\right) \ , \ n\geq 0}$ is a
martingale. Using the optimal sampling theorem we have

$$\E_x\sum\li_{k=0}^{\ntl}\al_k=\E_x\sum\li_{k=0}^{\ntl}
(\al_k-\bt_k+\bt_{k+1})+\E_x(\bt_0-\bt_{\ntl+1})
=\E_x(\bt_0-\bt_{\ntl+1}) \ .$$

The above argument shows that for the proof of $|(II)|\ra 0$
uniformly in $x\in [G]$ as $\ve \da 0$ it suffices to prove
$\play{\sup\li_{x\in \gm}\left|\sum\li_{n=0}^\infty
\E_x\al_n\right|}\ra 0$ uniformly in $x\in \gm$ as $\ve \da 0$.

Let us first show that $\E_{\nu}\al_0=0$. By Lemma 4.11 the Markov
chain $X_{\sm_n}^\ve$, $n \geq 0$ on $\gm$ is ergodic and has
invariant measure $\nu$. Therefore we have $\play{\lim\li_{n\ra
\infty}\dfrac{\sm_n}{n}=\E_{\nu}\sm_1}$. By ergodicity of the
process $X_t^\ve$ and self-adjointness of $L_0$ and $L_1$ we also
have
$$\lim\li_{t\ra \infty}\E_{\nu}\dfrac{1}{t}\play{\int_0^t A
f(H(X_s^\ve))ds=\int_{[G]}A f(H(x))dx \ .}$$ These two equalities
imply that $\play{\E_{\nu}\al_0=\E_{\nu}\int_0^{\sm_1}A
f(H(X_s^\ve))ds=(\E_{\nu}\sm_1)\cdot \int_{[G]}A f(H(x))dx}$.

We have
$$\int_{[G]}Af(H(x))dx=\text{Volume}(\cE)\cdot Af(O)+\sum\li_{k=1}^r \int_{I_k}M_k(H_k)\cL_k f(H_k)dH_k \ .$$
(The notations agree with those in the definition of the process
$Y_t$.)

Since $$\int_{I_k}M_k(H_k)\cL_k
f(H_k)dH_k=\int_{I_k}\dfrac{d}{dH_k}\left(M_k(H_k)\overline{\aone}(H_k)\dfrac{df}{dH_k}\right)dH_k=
\dfrac{1}{2}p_k\cdot \lim\li_{H_k\ra 0}\dfrac{df}{dH_k}(k,H_k) \ ,$$
we can use our boundary condition (2.1) to have
$\play{\int_{[G]}Af(H(x))dx=0}$ and therefore $\E_{\nu}\al_0=0$.

From the fact that $\E_{\nu}\al_0=0$ one first derives that
$\sup\li_{x\in \gm}\E_x \al_n$ decays to $0$ exponentially fast and
therefore $\sup\li_{x\in \gm}\left|\sum\li_{n=0}^\infty \E_x
\al_n\right|<\infty$. It also gives, for $x\in \gm$, that, for
$\sm^*\in \{\sm_1,\sm_2,...\}$ we have

$$\left|\sum\li_{n=0}^\infty \E_x \al_n\right|\leq \|Af\|_{\infty}\cdot\E_x\sm^*+
\text{Var}(m_{\sm^*}^x-\nu)\cdot \sup\li_{x\in
\gm}\left|\sum\li_{n=0}^\infty \E_x \al_n\right| \ .$$

Using Lemma 3.1 we see that for any $\dt>0$ we can choose $\sm^*$
such that

$$\sup\li_{x\in \gm}\left|\sum\li_{n=0}^\infty \E_x \al_n\right|\leq \|Af\|_{\infty}\cdot\dt+
\dt\cdot \sup\li_{x\in \gm}\left|\sum\li_{n=0}^\infty \E_x
\al_n\right| \ , $$ which proves that $\play{\sup\li_{x\in
\gm}\left|\sum\li_{n=0}^\infty \E_x\al_n\right|}\ra 0$ uniformly in
$x\in \gm$ as $\ve \da 0$. This implies that $|(II)|\ra 0$ uniformly
in $x\in [G]$ as $\ve \da 0$ and Theorem 2.1 follows. $\square$

\

\section{Auxiliary results needed in the proof of Theorem 2.1}

We establish in this section all the auxiliary results needed in
Section 3 for the proof of Theorem 2.1.

Let us make some further geometric constructions. Since we assumed
that the closure of all these $U_k$'s do not intersect each other we
see that for sufficiently small neighborhoods of these $U_k$'s they
also do not intersect each other. Without loss of generality let us
speak about one of these $U_k$'s. We remind that the matrix
$\azero(x)=(\azero_{ij}(x))_{1\leq i,j\leq d}$ is non-negative
definite inside $[G]$ and has rank $d$ on $[G]\backslash
\cup_{k=1}^r[U_k]$ and rank $(d-1)$ on $\cup_{k=1}^r
([U_k]\backslash \{x_k(m_k)\})$. At the points $x_k(m_k)$ the matrix
$\azero(x)$ is just $0$. From our assumptions we see that
$\azero(x)\grad H_k=0$ on $C_k(h)$ and $\e\cdot (\azero(x)\e)\geq
\underline{a}(x)|\e|^2_{\R^d}$ for any unit vector $\e$ tangent to
$C_k(h)$. Here $\underline{a}(x)>0$ for $x\in C_k(h)$ and $h\in
(m_k,0]$. The eigenvalue $\lb(x)=0$ for $\azero(x)$, $x\in
C_k(0)=\gm_k$ is simple and is the smallest one in the spectrum of
$\azero(x)$. For $x\in \gm_k$ the family of eigen-polynomials
$p(\lb;x)=\det(\lb I-a^{(0)}(x))$ pass through the origin. They are
transversal (i.e. not tangent) to the axis $p=0$. The transversality
is preserved under a small perturbation. From here one can see that
the eigenvalue $\lb(x)$ will remain simple and is still the smallest
one in the spectrum of all the matrices $\azero(x)$ as $x$ belongs
to a small neighborhood of $\gm_k$. We then see from implicit
function theorem that this eigenvalue $\lb(x)$ is a
$\contfunc^{(3)}$ function in a small enough neighborhood of
$\gm_k$. As a consequence, the unit eigenvector $\e_d(x)$ (for
different $k$ it is different vector fields but for simplicity of
notation we ignore that $k$ in our notation) corresponding to this
smallest eigenvalue is a $\contfunc^{(3)}$ vector field in a
neighborhood of $\gm_k$, with $\e_d(x)=\dfrac{\grad H_k (x)}{|\grad
H_k (x)|_{\R^d}}$ for $x\in \gm_k$. Let $X^x(t)$ be the integral
curve of this vector field. We let $\dfrac{d
X^x(t)}{dt}=\e_d(X^x(t))$, $X^x(0)=x\in \gm_k$. As we are working
within a small neighborhood of $\gm_k$ and $\e_d(x)$ in this
neighborhood is a $\contfunc^{(3)}$ vector field, being transversal
to $\gm_k$ when $x\in \gm_k$, we see that for $t\in [0,
\overline{h}]$ with $\overline{h}$ sufficiently small the points
$X^x(t)$ for fixed $t$ and all $x\in \gm_k$ form a surface
$\contfunc^{(3)}$ diffeomorphic to $\gm_k$. In this way we obtain an
extension of $H_k$ to a neighborhood of $U_k$ by letting
$H_k(X^x(t))=t$ for $t\in [0, \overline{h}]$. The Euclidean distance
from a point $X^x(t)$ to $\gm_k$ is $\geq \underline{d} \cdot t$ for
some $\underline{d}>0$. Let us denote by $C_k(+t)$ the level surface
$\{H_k=+t\}$ for $t\in [0, \overline{h}]$. Let
$\gmout_k=C_k(+\ve^{1/4})$ and $\gmoout_k=C_k(+2\ve^{1/4})$. Let
$\gmout=\cup_{k=1}^r \gmout_k$ and $\gmoout=\cup_{k=1}^r \gmoout_k$.
We can take $\ve$ small such that all $\gmoout_k$'s do not intersect
each other and do not touch $\pt G$. We denote by $\cE(\ve^{1/4})$
those points of $x\in \cE$ which lie outside the union of the
neighborhoods of the $U_k$'s bounded by $\gmout_k$, and we denote
$\cE(2\ve^{1/4})$ in a similar way.

We shall denote, for $x\in \cE(\ve^{1/4})$, the stopping time
$\sm(\ve^{1/4})$ to be the time when $X_t^\ve$ first hits $\gmout$.
Notice that by our assumption, for a point $x\in \cE(\ve^{1/4})$ we
have $$\sum\li_{i,j=1}^d \azero_{ij}(x)\xi_i\xi_j\geq \text{const}
\cdot \ve^{1/2} \sum\li_{i=1}^d \xi_i^2$$ for all
$(\xi_1,...,\xi_d)\in \R^d$ and some $\text{const}>0$.

\

Within the rest of this section implied positive constants denoted
by $C_i$'s will not be explicitly pointed out unless necessary.
Also, sometimes we use the same symbol $C$ to denote different
implied positive constants which are not important.

\

\textbf{Lemma 4.1.} \textit{For any $0<\varkappa<1/2$, for any $\ve$
small enough we have
$$\sup\li_{x\in [\cE(\ve^{1/4})]}
\E_x\sm(\ve^{1/4})\leq C\ve^{1/2-\varkappa}$$ for some $C>0$.}

\

\textbf{Proof.} Let $u^\ve(x,t)=\Prob_x(\sm(\ve^{1/4})>t)$. Then
$u^\ve(x,t)$ solves the problem

$$\left\{\begin{array}{l}
\dfrac{\pt u^\ve}{\pt t}=\left(\dfrac{1}{\ve}L_0+L_1\right)u^\ve \ ,
\\
u^\ve(y,0)=1 \text{ for } y \in \cE(\ve^{1/4}) \ ,
\\
u^\ve(y,t)=0 \text{ for } y \in \gmout \text{ and } t>0 \ ,
\\
\dfrac{\pt u^\ve}{\pt \gm^\ve}(y,t)=0 \text{ for } y \in \pt G \ .
\end{array}\right.$$

Our argument is a slight modification of the standard estimates of
heat kernel temporal decay (see \cite{[Fabes-Strook 1986]}). We
first multiply the equation that $u^\ve$ satisfies and we integrate
with respect to $x$ in $\cE(\ve^{1/4})$:

$$\begin{array}{l}
\play{\dfrac{d}{dt}\int_{\cE(\ve^{1/4})}(u^\ve)^2dx}
\\
\play{=\int_{\cE(\ve^{1/4})}u^\ve \grad \cdot
\left[\left(\dfrac{1}{\ve}a^{(0)}(x)+a^{(1)}(x)\right)\grad
u^\ve\right]dx}
\\
\play{=\int_{\cE(\ve^{1/4})} \grad \cdot \left[ u^\ve
\left(\dfrac{1}{\ve}a^{(0)}(x)+a^{(1)}(x)\right)\grad
u^\ve\right]dx-\int_{\cE(\ve^{1/4})}\left<\left(\dfrac{a^{(0)}(x)}{\ve}+a^{(1)}(x)\right)\grad
u^\ve, \grad u^\ve\right>_{\R^d}dx}
\\
\play{=-\int_{\cE(\ve^{1/4})}\left<\left(\dfrac{a^{(0)}(x)}{\ve}+a^{(1)}(x)\right)\grad
u^\ve, \grad u^\ve\right>_{\R^d}dx} \ .
\end{array}$$

The last step in the above calculation makes use of our boundary
condition:

$$\begin{array}{l}
\play{\int_{\cE(\ve^{1/4})} \grad \cdot \left[ u^\ve
\left(\dfrac{1}{\ve}a^{(0)}(x)+a^{(1)}(x)\right)\grad
u^\ve\right]dx}
\\
\play{=\int_{\pt G}u^\ve
\left(\dfrac{1}{\ve}a^{(0)}(x)+a^{(1)}(x)\right)\grad u^\ve\cdot
\textbf{n} dS+\int_{\gmout}u^\ve
\left(\dfrac{1}{\ve}a^{(0)}(x)+a^{(1)}(x)\right)\grad u^\ve\cdot
\textbf{n} dS}
\\
=0 \ .
\end{array}$$

Therefore we see that we have, for some $C>0$, that

$$\dfrac{d}{dt}\int_{\cE(\ve^{1/4})}(u^\ve)^2dx\leq
-\dfrac{C}{\ve^{1/2}}\int_{\cE(\ve^{1/4})}|\grad u^\ve|^2dx\ .$$

We can extend the function $u^\ve$ to the whole domain of $G$ so
that $u^\ve=0$ on $\cE\backslash \cE(\ve^{1/4})$. We can then apply
a variant of Poincar\'{e} inequality (see \cite[Lemma 1]{[Lieb et
al]}) so that we have

$$\int_{\cE(\ve^{1/4})}|\grad
u^\ve|^2dx\geq C\left(\int_{\cE(\ve^{1/4})}(u^\ve)^qdx\right)^{1/q}
$$ for $1<q<\dfrac{2d}{d-2}$. The constant $C>0$ in the above inequality can be chosen
independent of $\ve$.

Now we use the H\"{o}lder inequality so that, for
$\dfrac{1}{\al}+\dfrac{1}{\bt}=1$, $\al>0$, $\bt>0$ we obtain

$$\int_{\cE(\ve^{1/4})}(u^\ve)^2dx\leq
\left(\int_{\cE(\ve^{1/4})}u^\ve dx\right)^{1/\al}
\left(\int_{\cE(\ve^{1/4})}(u^\ve)^{(2-1/\al)\bt}dx\right)^{1/\bt} \
.$$

Combining the above two inequalities we see that we have, for
$\al>\dfrac{d+2}{4}$, that

$$\int_{\cE(\ve^{1/4})}(u^\ve)^2dx\leq
C\left(\int_{\cE(\ve^{1/4})}u^\ve dx\right)^{1/\al}
\left(\int_{\cE(\ve^{1/4})}|\grad u^\ve|^2dx\right)^{1-1/2\al} \ .$$

Thus we have

$$\dfrac{d}{dt}\int_{\cE(\ve^{1/4})}(u^\ve)^2dx\leq
-\dfrac{C}{\ve^{1/2}}
\dfrac{\play{\left(\int_{\cE(\ve^{1/4})}(u^\ve)^2dx\right)^{\frac{2\al}{2\al-1}}}}
{\play{\left(\int_{\cE(\ve^{1/4})}u^\ve
dx\right)^{\frac{2}{2\al-1}}}} \ .$$

Integrating in $x$ the equation that $u^\ve$ satisfies and making
use of the parabolic maximum principle it is possible to see that we
have

$$\int_{\cE(\ve^{1/4})}u^\ve dx\leq C$$ for some $C>0$.

This shows that we have

$$\dfrac{d}{dt}\int_{\cE(\ve^{1/4})}(u^\ve)^2dx\leq
-\dfrac{C}{\ve^{1/2}}
\play{\left(\int_{\cE(\ve^{1/4})}(u^\ve)^2dx\right)^{\frac{2\al}{2\al-1}}}
\ .$$

We shall denote $S(t)$ the semigroup generated by the operator
$\dfrac{1}{\ve}L_0+L_1$ with the prescribed boundary conditions as
in the problem for $u^\ve$. From co-normal condition one can check
that $S(t)$ is self-adjoint. From the above inequality it is
standard to deduce
$$\|S(t)\|_{L^1(\cE(\ve^{1/4}))\ra L^2(\cE(\ve^{1/4}))}\leq \dfrac{C}{(\ve^{-1/2} t)^{q}}$$ for any
$q>d/4$. By self-adjointness of $S(t)$ we see that
$$\|S(t)\|_{L^2(\cE(\ve^{1/4}))\ra L^\infty(\cE(\ve^{1/4}))}\leq
\dfrac{C}{(\ve^{-1/2} t)^{q}} \ .$$ Since $S(t)=S(t/2)\circ S(t/2)$
we have $$\|S(t)\|_{L^1(\cE(\ve^{1/4}))\ra
L^\infty(\cE(\ve^{1/4}))}\leq \dfrac{C}{(\ve^{-1} t^2)^{q}}$$ for
any $q>d/4$. In particular, this means that we have
$$\sup\li_{x\in [\cE(\ve^{1/4})]}\Prob_x(\sm(\ve^{1/4})>\ve^{1/2-\kp})\leq C\ve^{2\kp q}$$
for $q>d/4$.

By strong Markov property of the process $X_t^\ve$ we see that
$$\begin{array}{l}
\E_x\sm(\ve^{1/4})
\\
=\play{\int_0^\infty \Prob_x(\sm(\ve^{1/4})>t)dt}
\\
\leq \ve^{1/2-\varkappa} \sum\li_{n=0}^\infty
\Prob_x(\sm(\ve^{1/4})>n\ve^{1/2-\varkappa})
\\
\leq \ve^{1/2-\varkappa} \sum\li_{n=0}^\infty (\sup\li_{x\in
[\cE(\ve^{1/4})]}\Prob_x(\sm(\ve^{1/4})>\ve^{1/2-\varkappa}))^n
\\
= \dfrac{\ve^{1/2-\varkappa}}{1-\sup\li_{x\in
[\cE(\ve^{1/4})]}\Prob_x(\sm(\ve^{1/4})>\ve^{1/2-\varkappa})} \leq
C\ve^{1/2-\varkappa}
\end{array}$$ for $\ve$ small enough. This implies the statement of the Lemma. $\square$

\

We shall denote by $\cS_k([0,\ve^{1/4}])$ the closed set bounded by
the surfaces $\gmout_k$ and $\gm_k$ and by
$\cS([0,\ve^{1/4}])=\cup_{k=1}^r\cS_k([0,\ve^{1/4}])$. We denote
$\cS_k([0,2\ve^{1/4}])$ and $\cS([0,2\ve^{1/4}])$ in a similar way
by replacing $\gmout_k$ by $\gmoout_k$.

Following the geometric construction stated before Lemma 4.1, for
$\ve>0$ small enough, and each $k=1,...,r$, at any point $x\in
\cS_k([0,2\ve^{1/4}])$ one can find $\e_d(x)=\dfrac{\grad
H_k(x)}{|\grad H_k(x)|_{\R^d}}$ and for any unite vector $\e$ such
that $\e\cdot \e_d(x)=0$ we have $\e\cdot (\azero(x)\e)\geq
\underline{a}|\e|_{\R^d}^2=\underline{a}$ for some
$\underline{a}>0$. Also $\azero(x)\e_d(x)=\lb(x)\e_d(x)$. The
eigenvalue $\lb(x)$ is in $\contfunc^{(3)}(\cS_k([0,2\ve^{1/4}]))$
with $\lb|_{\gm_k}=0$ and $\lb(x)>0$ for $x\in
\cS_k([0,2\ve^{1/4}])\backslash \gm_k$. Furthermore, for $\ve$ small
enough we have $C_1\cdot\text{dist}^2(x,\gm_k) \leq \lb(x)\leq
C_2\cdot\text{dist}^2(x,\gm_k)$ for some $C_1,C_2>0$ and $x\in
\cS([0,2\ve^{1/4}])$.

For each $x\in \cS_k([0, 2\ve^{1/4}])$, we can find a point
$\vphi^k(x)\in \gm_k$ such that $X^{\vphi^k(x)}(H_k(x))=x$ for the
flow $X^x(t)$ introduced in the geometric construction before Lemma
4.1. Let us introduce a local coordinate
$(\vphi^k_1,...,\vphi^k_{d-1}, H_k)$ in a small neighborhood of a
point $x\in \cS_k([0,2\ve^{1/4}])$. We take $H_k=H_k(x)$, which is
the extended first integral of $H_k$ to $\cS_k([0,2\ve^{1/4}])$; and
we take
$(\vphi^k_1,...,\vphi^k_{d-1})=(\vphi^k_1(x),...,\vphi^k_{d-1}(x))$
to be the local coordinate for the point $\vphi^k(x)$ on $\gm_k$. In
the more or less simpler case we can arrange the coordinate
$(\vphi^k_1(x),...,\vphi^k_{d-1}(x), H_k(x))$ in such a way that
this new coordinate system is an orthogonal curvilinear coordinate
system. (We will discuss the general case a bit later.) The metric
tensor of this new coordinate system is given by
$ds^2=E_1(x)(d\vphi^k_1(x))^2+...+E_{d-1}(x)(d\vphi^k_{d-1}(x))^2+E_d(x)(dH_k(x))^2$.
Here the functions $0<C_3<E_1(x),...,E_d(x)<C_4<\infty$ are in class
$\contfunc^{(3)}(\cS_k([0,2\ve^{1/4}]))$ with bounded derivatives.
We notice that by our geometric construction we have $C_5\cdot
H_k^2(x) \leq \lb(x)\leq C_6\cdot H_k^2(x)$. We shall let $\e_i(x)$
be the unit tangent vector on the axis curve corresponding to
$\vphi_i^k$ for $1\leq i\leq d-1$; $\e_d(x)$ be the unit tangent
vector on the axis curve corresponding to $H_k$. The vectors
$\e_1(x),...,\e_d(x)$ are basis vectors.

The theory of orthogonal curvilinear coordinate system (see, for
example, \cite[Ch.14]{[Zorich]}) tells us that for a differentiable
function $f$ on $\cS_k([0,2\ve^{1/4}])$ we have

$$\grad f (x)=\sum\li_{i=1}^{d-1}\dfrac{1}{\sqrt{E_i(x)}}\dfrac{\pt f}{\pt \vphi^k_i}(x)\e_i(x)
+\dfrac{1}{\sqrt{E_d(x)}}\dfrac{\pt f}{\pt H_k}(x) \e_d(x) \ ,$$ and
for a differentiable vector field $\mathbf{B}(x)=\sum\li_{i=1}^d
B^i(x)\e_i(x)$ on $\cS_k([0,2\ve^{1/4}])$ we have

$$\grad \cdot \mathbf{B}(x)= \dfrac{1}{\sqrt{\prod_{i=1}^d E_i(x)}}\left[\sum\li_{i=1}^{d-1}
\dfrac{\pt}{\pt \vphi_i^k}\left(\sqrt{\dfrac{\prod_{j=1}^d
E_j(x)}{E_i(x)}}B^i(x)\right)+\dfrac{\pt}{\pt
H_k}\left(\sqrt{\dfrac{\prod_{j=1}^d
E_j(x)}{E_d(x)}}B^d(x)\right)\right] \ .$$

Consider a function (so called "barrier function", see
\cite{[Hasminskii]} and \cite[Ch.3]{[F red book]}) $u_k(x)\in
\contfunc^{(2)}(\cS_k([0,2\ve^{1/4}]))$ which depends only on $H_k$
and is a constant on each level surface $\{H_k=\text{const}\}$. We
can write $u_k(x)=u_k(H_k)$ and we apply the above two formulas to
get

$$\begin{array}{l}
\left(\dfrac{1}{\ve}L_0+L_1\right)u_k(x)
\\
=\left(\dfrac{1}{2\ve}\grad \cdot (\azero(x)\grad
u_k(x))+\dfrac{1}{2}\grad \cdot (\aone(x)\grad u_k(x))\right)
\\
=\dfrac{1}{\sqrt{\prod_{i=1}^d
E_i(x)}}\left[\dfrac{1}{2}\dfrac{\pt}{\pt
H_k}\left(\sqrt{\dfrac{\prod_{i=1}^d
E_i(x)}{E_d^2(x)}}\left(\dfrac{\lb(x)}{\ve}+\mu_d(x)\right)\dfrac{du_k}{dH_k}(H_k)\right)\right.
\\
\left. \ \ \ \ \ +\sum\li_{i=1}^{d-1}\dfrac{1}{2}\dfrac{\pt}{\pt
\vphi_i^k}\left(\sqrt{\dfrac{\prod_{i=1}^d
E_i(x)}{E_d(x)E_i(x)}}\mu_i(x)\right)\cdot
\dfrac{du_k}{dH_k}(H_k)\right] \ .
\end{array}$$

Here the functions $\mu_1(x),...,\mu_d(x)$ are defined via the
relation $\aone(x) \e_d(x)=\mu_1(x)\e_1(x)+...+\mu_d(x)\e_d(x)$.
These functions are in $\contfunc^{(3)}(\cS([0,2\ve^{1/4}]))$ with
bounded derivatives. Notice that since $L_1$ is strictly elliptic,
the matrix $\aone(x)$ is positive definite, and therefore the
function $\mu_d(x)$ is uniformly bounded from below by a certain
positive constant.

For simplicity of notation let us define $A(x)=\sqrt{\prod_{i=1}^d
E_i(x)}$ and $A_i(x)=\dfrac{A(x)}{\sqrt{E_i(x)E_d(x)}}$ for
$i=1,...,d$. These functions are strictly positive (with uniform
lower bound) in $\contfunc^{(2)}(\cS([0,2\ve^{1/4}]))$ with bounded
derivatives. Under this notation we can write

$$\begin{array}{l}
\left(\dfrac{1}{\ve}L_0+L_1\right)u_k(x)
\\
=\dfrac{1}{A(x)}\left[\dfrac{1}{2}\dfrac{\pt}{\pt
H_k}\left(A_d(x)\left(\dfrac{\lb(x)}{\ve}+\mu_d(x)\right)\dfrac{du_k}{dH_k}(H_k)\right)
+\sum\li_{i=1}^{d-1}\dfrac{1}{2}\dfrac{\pt}{\pt
\vphi_i^k}\left(A_i(x)\mu_i(x)\right)\cdot
\dfrac{du_k}{dH_k}(H_k)\right] \ .
\end{array}$$

As a further simplification we shall define
$$\dfrac{1}{2}A_d(x)\lb(x)=K_1(x) \ ,$$
$$\dfrac{1}{2}A_d(x)\mu_d(x)=K_2(x) \ ,$$
$$\sum\li_{i=1}^{d-1}\dfrac{1}{2}\dfrac{\pt}{\pt
\vphi_i^k}(A_i(x)\mu_i(x))=K_3(x) \ .$$ We have
$$\play{\left(\dfrac{1}{\ve}L_0+L_1\right)u_k(x)=\dfrac{1}{A(x)}\left[\dfrac{\pt}{\pt
H_k}\left(\left(\dfrac{K_1(x)}{\ve}+K_2(x)\right)\dfrac{du_k}{dH_k}(H_k)\right)+K_3(x)\dfrac{du_k}{dH_k}(H_k)\right]}
\ . \eqno(4.1)$$

For a point $x\in \cS_k([0,2\ve^{1/4}])$ and $\ve$ small enough we
have $$C_7 H_k^2(x)\leq K_1(x)\leq C_8 H_k^2(x) \ ; \eqno(4.2)$$
$$C_9 H_k(x)\leq \dfrac{\pt}{\pt H_k}(K_1(x))\leq C_{10}H_k(x) \ ; \eqno(4.3)$$
$$0<C_{11}<K_2(x)<C_{12}<\infty \ ; \eqno(4.4)$$ $$\left|\dfrac{\pt}{\pt
H_k}(K_2(x))\right|\leq C_{13}<\infty  \ ; \eqno(4.5)$$
$$|K_3(x)|\leq C_{14}<\infty \ . \eqno(4.6)$$

We also notice that since we are working in $\cS_k([0,2\ve^{1/4}])$
and $\ve$ is small, the functions $A_d(x)=A_d(\vphi^k,H_k)$ and
$\lb(x)=\lb(\vphi^k,H_k)$ have Taylor expansions
$$A_d(\vphi^k,H_k)=A_d(\vphi^k,0)+O(H_k)\ ,$$
$$\lb(\vphi^k,H_k)=\dfrac{1}{2}\dfrac{\pt^2 \lb}{\pt H_k^2}(\vphi^k,0)H_k^2+O(H_k^3)\ .$$

Therefore we see that for $x\in \cS_k([0,2\ve^{1/4}])$ we have

$$K_1(x)=C_k(\vphi^k)H_k^2+O(H_k^3) \eqno(4.7)$$
with a certain positive function $C_k(\vphi^k)$.

\

In the general case the axis curve corresponding to $H_k$ will be
orthogonal to those corresponding to the $\vphi_i^k$'s, but the axis
curves corresponding to the $\vphi_i^k$'s are not necessarily
orthogonal. The calculation will be more bulky since the metric
tensor have cross terms with respect to the coordinate
$\vphi_i^k$'s, but the essence is the same as it is only important
to have the axis curves corresponding to $H_k$ being orthogonal to
those corresponding to the $\vphi_i^k$'s. To be more precise, let
$(g_{ij})_{1\leq i,j\leq d}$ be the metric tensor corresponding to
the (local) coordinate system $(\vphi_1^k,...,\vphi_{d-1}^k, H_k)$.
Let $\e_i(x)$ be the unit tangent vector on the axis curve
corresponding to $\vphi_i^k$ for $1\leq i\leq d-1$; $\e_d(x)$ be the
unit tangent vector on the axis curve corresponding to $H_k$. The
vectors $\e_1(x),...,\e_d(x)$ are basis vectors. We have
$g_{id}=g_{di}=0$ for $i=1,...,d-1$ and $g_{dd}>0$. Let
$(g^{ij})_{1\leq i,j \leq d}$ be the dual tensor, i.e.,
$(g^{ij})_{1\leq i,j\leq d}$ is the inverse matrix of
$(g_{ij})_{1\leq i,j\leq d}$. We have $g^{id}=g^{di}=0$ for $1\leq
i\leq d-1$ and $g^{dd}=\dfrac{1}{g_{dd}}$. For $u_k=u_k(H_k)$ we
have

$$\grad u_k(x)=\dfrac{1}{\sqrt{g_{dd}(x)}}\dfrac{du_k}{dH_k}\e_d(x) \ ,$$
and
$$\azero(x)\grad u_k(x)=\dfrac{\lb(x)}{\sqrt{g_{dd}(x)}}\dfrac{du_k}{dH_k}\e_d(x) \ ;$$
$$\aone(x)\grad u_k(x)=\dfrac{\mu_d(x)}{\sqrt{g_{dd}(x)}}\dfrac{du_k}{dH_k}\e_d(x)+
\dfrac{1}{\sqrt{g_{dd}(x)}}\dfrac{du_k}{dH_k}(\mu_1(x)\e_1(x)+...+\mu_{d-1}(x)\e_{d-1}(x))
\ .$$

Here, as before, we have
$\aone(x)\e_d(x)=\mu_1(x)\e_1(x)+...+\mu_d(x)\e_d(x)$. We shall then
apply a general formula that for a vector field
$\mathbf{B}(x)=\sum\li_{i=1}^dB^i(x)\e^i(x)$ we have $$\grad \cdot
\mathbf{B}(x)=\dfrac{1}{\sqrt{g(x)}}\sum\li_{i=1}^{d-1}\dfrac{\pt}{\pt
\vphi_i^k}(B^i(x)\sqrt{g^{ii}(x)}\sqrt{g(x)})+\dfrac{1}{\sqrt{g(x)}}\dfrac{\pt}{\pt
H_k}(B^d(x)\sqrt{g^{dd}(x)}\sqrt{g(x)}) \ .$$ Here
$g(x)=\det(g_{ij}(x))$. The basis $\e^1(x),...,\e^d(x)$ is the
reciprocal basis (normalized) dual to $\e_1(x),...,\e_d(x)$, i.e.,
$(\e_i,\e^j)_{(g_{ij})}=\dt_{ij}$ with respect to the inner product
$(\bullet,\bullet)_{(g_{ij})}$ defined by the metric tensor
$(g_{ij})$. By the fact that the metric tensor has no cross terms
between $H_k$ and $\vphi_i^k$'s, we actually have $\e^d(x)=\e_d(x)$
and
$\text{span}\{\e_1(x),...,\e_{d-1}(x)\}=\text{span}\{\e^1(x),...,\e^{d-1}(x)\}$.

We then see that the operator $\dfrac{1}{\ve}L_0+L_1$ applied to
$u_k(x)=u_k(H_k)$ will result in a formula which is the same as
(4.1). The functions $K_1(x)$, $K_2(x)$ and $K_3(x)$ will somehow be
different but they still satisfy the conditions (4.2)\,--\,(4.7).

\

Let $\zt([0,2\ve^{1/4}])$ be the first time when the process
$X_t^\ve$, starting from a point $x\in \cS([0,2\ve^{1/4}])$, hits
$\gm$ or $\gmoout$.

\

\textbf{Lemma 4.2.} \textit{We have }
$$\sup\li_{x\in
\cS([0,2\ve^{1/4}])} \E_x\zt([0,2\ve^{1/4}])\leq \ C\ve^{3/4}
$$ \textit{for some $C>0$.}

\

\textbf{Proof.} Let $$K_4(H_k)=\dfrac{H_k^2}{\ve}+1 \ $$ and
$$K_5(x)=\dfrac{1}{K_4(H_k(x))}\left(\dfrac{K_1(x)}{\ve}+K_2(x)\right) \ .$$

By (4.2) and (4.4) we can estimate
$$C_{15}\leq K_5(x)\leq
C_{16} \eqno(4.8)$$ for $C_{15}, C_{16}>0$, $x\in \cS_k([0,
2\ve^{1/4}])$ and $\ve$ small enough.

By (4.7) we see that

$$\begin{array}{l}
K_5(x)
\\
\play{=\dfrac{K_1(x)+\ve K_2(x)}{H_k^2+\ve}}
\\
\play{=\dfrac{C_k(\vphi^k(x))H_k^2+O(H_k^3)+\ve K_2(x)}{H_k^2+\ve}}
\\
\play{=C_k(\vphi^k(x))+H_k\cdot\dfrac{O(H_k^2)+\dfrac{\ve}{H_k}
(K_2-C_k(\vphi^k(x)))}{H_k^2+\ve}} \ .
\end{array}$$

Since $0\leq H_k\leq 2\ve^{1/4}$ we see that as $\ve$ is small we
have

$$\left|\dfrac{O(H_k^2)+\dfrac{\ve}{H_k}
(K_2-C_k(\vphi^k(x)))}{H_k^2+\ve}\right|\leq C_{17} \ .$$

Thus as $0\leq H_k\leq 2\ve^{1/4}$ we see that
$$\left|\dfrac{\pt}{\pt H_k}(K_5(x))\right|\leq C_{17} \ .$$

Let the barrier function $u_k^{(1)}(x)=u_k^{(1)}(H_k)$ be defined by

$$u_k^{(1)}(H_k)=\int_0^{H_k}\dfrac{K_6(\ve)-y}{K_4(y)}dy$$ with
$$K_6(\ve)=\left(\int_0^{2\ve^{1/4}}\dfrac{dy}{K_4(y)}\right)^{-1}
\left(\int_0^{2\ve^{1/4}}\dfrac{y dy}{K_4(y)}\right) \ .$$

It is easy to check that $$u_k^{(1)}(0)=u_k^{(1)}(2\ve^{1/4})=0 \
.$$

We can estimate $K_6(\ve)\leq 2\ve^{1/4}$ and we have

$$\play{\int_0^{H_k}\dfrac{dy}{K_4(y)}}=\int_{0}^{H_k}\dfrac{dy}{\dfrac{y^2}{\ve}+1}
=\ve^{1/2}\arctan(H_k\ve^{-1/2})\leq C_{18}\ve^{1/2} \ .
\eqno(4.9)$$

This gives the estimates $$0\leq u_k^{(1)}(H_k)\leq C_{19}\ve^{3/4}
\eqno(4.10)$$ and $$\left|\dfrac{du_k^{(1)}}{dH_k}(H_k)\right|\leq
C_{20}\ve^{1/4} \eqno(4.11)$$ for $0\leq H_k\leq 2\ve^{1/4}$. Apply
(4.1) to the function $u_k^{(1)}$ we can see, using (4.11), that,

$$\begin{array}{l}
\play{\left(\dfrac{1}{\ve}L_0+L_1\right)u_k^{(1)}(x)}
\\
\play{=\dfrac{1}{A(x)}\left[\dfrac{\pt}{\pt
H_k}\left(K_5(x)K_4(H_k(x))\dfrac{du_k^{(1)}}{dH_k}(H_k(x))\right)+K_3(x)\dfrac{du_k^{(1)}}{dH_k}(H_k(x))\right]}
\\
\play{\leq\dfrac{1}{A(x)}\left[\dfrac{\pt}{\pt
H_k}\left((K_6(\ve)-H_k(x))K_5(x)\right)+C_{21}\ve^{1/4}\right]}
\\
\play{=\dfrac{1}{A(x)}\left[-K_5(x)+\dfrac{\pt}{\pt H_k}(
K_5(x))(K_6(\ve)-H_k(x))+C_{21}\ve^{1/4}\right]\leq -C_{22}}
\end{array} \eqno(4.12)$$
for $x\in \cS([0,2\ve^{1/4}])$ and $\ve$ small enough.

We notice that this process $X_t^\ve$ before hitting $\gm$ or
$\gmoout$ is restricted to one of the $\cS_k([0,2\ve^{1/4}])$'s and
the bound (4.12) can be made uniform in $k$.

Now we apply It\^{o}'s formula to the function $u_k^{(1)}$
constructed above up to the stopping time $\zt([0,2\ve^{1/4}])$.
Taking expectation we get
$$u_k^{(1)}(x)=-\E_x\play{\int_0^{\zt([0,2\ve^{1/4}])}\left(\dfrac{1}{\ve}L_0+L_1\right)u_k^{(1)}
(X_s^\ve) ds}\geq C_{22}\E_x\zt([0,2\ve^{1/4}]) \ .\eqno(4.13) $$
From (4.10) and (4.13) we see that the statement of this Lemma
follows. $\square$

\

\textbf{Lemma 4.3.} \textit{For $x\in \gmout$ we have}
$$\Prob_x(X_{\zt([0,2\ve^{1/4}])}^\ve\in \gm)\geq C\ve^{1/4}$$ \textit{for some $C>0$.}

\

\textbf{Proof.} Let

$$K_7(H_k)=\max\li_{x\in \cS_k([0,2\ve^{1/4}]),
H_k(x)=H_k}\dfrac{\dfrac{\pt}{\pt
H_k}\left(\dfrac{K_1(x)}{\ve}+K_2(x)\right)+K_3(x)}{\dfrac{K_1(x)}{\ve}+K_2(x)}
\ .$$

Let, for a fixed $H_k\in[0,2\ve^{1/4}]$, the above maximum be
achieved at a point $(\vphi^k,H_k)$.

We will prove in Lemma 4.4 that we have the auxiliary estimate

$$\left|K_7(H_k)-\dfrac{2C_k(\vphi^k)H_k}{C_k(\vphi^k)H_k^2+\ve K_2(\vphi^k, H_k)}\right|\leq C_{23} \ . \eqno(4.14)$$

Let the barrier function $u_k^{(2)}(x)=u_k^{(2)}(H_k)$ be defined by

$$u_k^{(2)}(H_k)=1-\dfrac{\play{\int_0^{H_k}\exp\left(-\int_0^y
K_7(z)dz\right)dy}}{\play{\int_0^{2\ve^{1/4}}\exp\left(-\int_0^y
K_7(z)dz\right)dy}} \ .$$

It is easy to see that we have
$$u_k^{(2)}(0)=1 \ , \ u_k^{(2)}(2\ve^{1/4})=0 \ .$$

Apply formula (4.1) we can see that

$$\begin{array}{l}
\play{\left(\dfrac{1}{\ve}L_0+L_1\right)u_k^{(2)}(x)}
\\
\play{=\dfrac{1}{A(x)}\left[\left(\dfrac{K_1(x)}{\ve}+K_2(x)\right)\dfrac{d^2u_k^{(2)}}{dH_k^2}(H_k)
+\left(\dfrac{\pt}{\pt
H_k}\left(\dfrac{K_1(x)}{\ve}+K_2(x)\right)+K_3(x)\right)\dfrac{du_k^{(2)}}{dH_k}(H_k)\right]}
\\
\geq 0 \ .
\end{array} \eqno(4.15)$$

However, by (4.14), (4.4) and the property of $C_k(\vphi^k)$ in
(4.7) we can estimate

$$\int_{H_k}^{2\ve^{1/4}}\exp\left(-\int_0^y K_7(z)dz\right)dy\geq C_{24}
\int_{H_k}^{2\ve^{1/4}}\exp\left(-\int_0^y
\dfrac{2z}{z^2+C_{25}\ve}dz\right)dy \ , \eqno(4.16)$$

$$\int_0^{2\ve^{1/4}}\exp\left(-\int_0^y K_7(z)dz\right)dy\leq C_{26}
\int_0^{2\ve^{1/4}}\exp\left(-\int_0^y
\dfrac{2z}{z^2+C_{27}\ve}dz\right)dy \ . \eqno(4.17)$$

By (4.16) and (4.17), and taking into account the auxiliary
estimates made in Lemma 4.5 we see that
$$u_k^{(2)}(\ve^{1/4})\geq C_{28}\ve^{1/4} \ . \eqno(4.18)$$

This bound (4.18) can actually be made uniform in $k$. We can apply
It\^{o}'s formula to the function $u_k^{(2)}$ constructed above up
to the stopping time $\zt([0,2\ve^{1/4}])$. Taking expectation we
get
$$\Prob_x(X_{\zt([0,2\ve^{1/4}])}^\ve\in
\gm)-u_k^{(2)}(\ve^{1/4})=\E_x\play{\int_0^{\zt([0,2\ve^{1/4}])}\left(\dfrac{1}{\ve}L_0+L_1\right)u_k^{(2)}
(X_s^\ve) ds}\geq 0 \  \eqno(4.19)$$ for $x\in \gmout$. Now (4.18)
and (4.19) imply the statement of this Lemma. $\square$

\

\textbf{Lemma 4.4.} \textit{For a fixed $H_k\in [0,2\ve^{1/4}]$ and
the corresponding $\vphi^k$ defined as in the proof of Lemma 4.3, we
have}

$$\left|K_7(H_k)-\dfrac{2C_k(\vphi^k)H_k}{C_k(\vphi^k)H_k^2+\ve K_2(\vphi^k, H_k)}\right|\leq C$$
\textit{for some $C>0$.}

\

\textbf{Proof.} Using (4.7), we can write

$$K_7(H_k)=\dfrac{2C_k(\vphi^k)H_k+O(H_k^2)+
\ve K_{2,3}(\vphi^k,H_k)}{C_k(\vphi^k)H_k^2+O(H_k^3)+\ve
K_2(\vphi^k,H_k)} \ .$$

Here $K_{2,3}(\vphi^k,H_k)$ is a bounded function. We then have

$$\begin{array}{l}
\play{\left|K_7(H_k)-\dfrac{2C_k(\vphi^k)H_k}{C_k(\vphi^k)H_k^2+\ve
K_2(\vphi^k, H_k)}\right|}
\\
= \play{\left|\dfrac{2C_k(\vphi^k)H_k+O(H_k^2)+ \ve
K_{2,3}(\vphi^k,H_k)}{C_k(\vphi^k)H_k^2+O(H_k^3)+\ve
K_2(\vphi^k,H_k)}-\dfrac{2C_k(\vphi^k)H_k}{C_k(\vphi^k)H_k^2+\ve
K_2(\vphi^k, H_k)}\right|}
\\
\leq \play{\left|\dfrac{O(H_k^2)+\ve
K_{2,3}(\vphi^k,H_k)}{C_k(\vphi^k)H_k^2+O(H_k^3)+\ve
K_2(\vphi^k,H_k)}\right|+}
\\
\play{\ \ \ \ \ \ \
+\left|\dfrac{2C_k(\vphi^k)H_k}{C_k(\vphi^k)H_k^2+O(H_k^3)+\ve
K_2(\vphi^k,H_k)}-\dfrac{2C_k(\vphi^k)H_k}{C_k(\vphi^k)H_k^2+\ve
K_2(\vphi^k,H_k)}\right|}
\\
= \play{\left|\dfrac{O(H_k^2)+\ve
K_{2,3}(\vphi^k,H_k)}{C_k(\vphi^k)H_k^2+O(H_k^3)+\ve
K_2(\vphi^k,H_k)}\right|+}
\\
\play{\ \ \ \ \ \ \
+\left|\dfrac{2C_k(\vphi^k)H_k}{C_k(\vphi^k)H_k^2+\ve
K_2(\vphi^k,H_k)}\cdot
\dfrac{O(H_k^3)}{C_k(\vphi^k)H_k^2+O(H_k^3)+\ve
K_2(\vphi^k,H_k)}\right|\leq C \ .}
\end{array}$$ $\square$

\

\textbf{Lemma 4.5.} \textit{We have}

$$\int_{\ve^{1/4}}^{2\ve^{1/4}}\exp\left(-\int_0^y \dfrac{2z}{z^2+C\ve}dz\right)dy\geq C_{29}\ve^{3/4} \ ,$$

$$C_{31}\ve^{1/2}\geq \int_{0}^{C_{32}\ve^{1/4}}\exp\left(-\int_0^y \dfrac{2z}{z^2+C\ve}dz\right)dy\geq C_{30}\ve^{1/2} \ .$$

\

\textbf{Proof.} Evaluating the integrals, we have

$$\int_0^y \dfrac{2z}{z^2+C\ve}dz=\ln\left(\dfrac{y^2+C\ve}{C\ve}\right) \ ,$$

$$\int_a^b \exp\left(-\int_0^y\dfrac{2z}{z^2+C\ve}dz\right)dy=
\sqrt{C}\ve^{1/2}\left(\arctan(\dfrac{b}{\sqrt{C}\ve^{1/2}})-\arctan(\dfrac{a}{\sqrt{C}\ve^{1/2}})\right)
\ .$$

If $a=0$ and $b=C_{32}\ve^{1/4}$ we already get the second
inequality of this Lemma. Now suppose $a=\ve^{1/4}$ and
$b=2\ve^{1/4}$. We shall make use of an asymptotic expansion of
$\arctan(y)$ as $y\ra \infty$:

$$\arctan (y)=\dfrac{\pi}{2}-\dfrac{1}{y}+O(\dfrac{1}{y^2}) \ \text{ as } \ y \ra \infty \ .$$

This gives

$$\begin{array}{l}
\play{\int_{\ve^{1/4}}^{2\ve^{1/4}}\exp\left(-\int_0^y
\dfrac{2z}{z^2+C\ve}dz\right)dy}
\\
\geq
C_{33}\ve^{1/2}\left(-\dfrac{1}{2\ve^{-1/4}}+\dfrac{1}{\ve^{-1/4}}+O(\ve^{1/2})\right)\geq
C_{29}\ve^{3/4} \ .
\end{array}$$ $\square$

\

\textbf{Lemma 4.6.} \textit{We have}

$$\lim\li_{\ve\da 0}\sup\li_{x\in \gmout}\E_x\sm=0 \ .$$

\

\textbf{Proof.} Lemmas 4.1, 4.2 and 4.3 imply the statement of this
Lemma. For $x\in \gmout$ we have

$$\begin{array}{l}
\E_x\sm
\\
=\E_x\zt([0,2\ve^{1/4}])\1(X_{\zt([0,2\ve^{1/4}])}^\ve\in
\gm)+\E_x(\zt([0,2\ve^{1/4}])+\E_{X_{\zt([0,2\ve^{1/4}])}^\ve}\sm)\1(X_{\zt([0,2\ve^{1/4}])}^\ve\in
\gmoout)
\\
\leq \sup\li_{x\in
\cS([0,2\ve^{1/4}])}\E_x\zt([0,2\ve^{1/4}])+\E_x(\sup\li_{y\in
\gmoout}\E_y \sm(\ve^{1/4})+\sup\li_{x\in \gmout}\E_x
\sm)\1(X_{\zt([0,2\ve^{1/4}])}^\ve\in \gmoout)
\\
\leq \sup\li_{x\in
\cS([0,2\ve^{1/4}])}\E_x\zt([0,2\ve^{1/4}])+(\sup\li_{y\in
\gmoout}\E_y \sm(\ve^{1/4})+\sup\li_{x\in \gmout}\E_x
\sm)\Prob_x(X_{\zt([0,2\ve^{1/4}])}^\ve\in \gmoout) \ .
\end{array} \eqno(4.20) $$

Taking a sup over all $x\in \gmout$ we get

$$\sup\li_{x\in\gmout}\E_x\sm\leq \dfrac{\sup\li_{x\in
\cS([0,2\ve^{1/4}])}\E_x\zt([0,2\ve^{1/4}])+\sup\li_{y\in
\gmoout}\E_y \sm(\ve^{1/4})\cdot
\Prob_x(X_{\zt([0,2\ve^{1/4}])}^\ve\in
\gmoout)}{\Prob_x(X_{\zt([0,2\ve^{1/4}])}^\ve\in \gm)}\ . $$

Using Lemmas 4.1, 4.2 and 4.3 we see that the statement of this
Lemma follows. (We choose $\varkappa=1/8$ in Lemma 4.1.) $\square$

\

\textbf{Lemma 4.7.} \textit{We have}

$$\lim\li_{\ve\da 0}\sup\li_{x\in \cS([0,\ve^{1/4}])}\E_x\sm=0 \ .$$

\

\textbf{Proof.} This is a consequence of Lemmas 4.1, 4.2 and 4.6.
$\square$

\

\textbf{Lemma 4.8.} \textit{We have}

$$\lim\li_{\ve\da 0}\sup\li_{x\in [\cE]}\E_x\sm=0 \ .$$

\

\textbf{Proof.} This is a consequence of Lemmas 4.1 and 4.7.
$\square$

\

We shall denote by $\cS_k([-\ve^{1/2},\ve^{1/4}])$ the closed set
bounded by the surfaces $\gmin_k$ and $\gmout_k$ and by
$\cS([-\ve^{1/2},\ve^{1/4}])=\cup_{k=1}^r\cS_k([-\ve^{1/2},\ve^{1/4}])$.
We notice that by the same reason as before, a coordinate
$(\vphi_1^k,...,\vphi_{d-1}^k, H_k)$ exists in
$\cS_k([-\ve^{1/2},\ve^{1/4}])$. We denote $\cS_k([0,\ve^{1/4}])$,
$\cS([0,\ve^{1/4}])$ (replacing $\gmin_k$ by $\gm_k$) and
$\cS_k([-\ve^{1/2},0])$, $\cS([-\ve^{1/2},0])$ (replacing $\gmout_k$
by $\gm_k$) in a similar way.

\

\textbf{Lemma 4.9.} \textit{We have $$\lim\li_{\ve\da
0}\sup\li_{x\in \gm} \E_x\tau=0 \ .$$}

\

\textbf{Proof.} The proof of this lemma is very similar to and is a
bit simpler than that of Lemma 4.6. We shall construct two barrier
functions $u_k^{(3)}$ (for estimating the exit time from
$\cS_k([-\ve^{-1/2},\ve^{1/4}])$) and $u_k^{(4)}$ (for the
probability of hitting $\gmin$).

For the construction of $u_k^{(3)}$ all the arguments of Lemma 4.2
can be carried here with $\gm$ replaced by $\gmin$, $\gmout$
replaced by $\gm$ and $\gmoout$ replaced by $\gmout$. We are working
now with $\cS_k([-\ve^{1/2},\ve^{1/4}])$ and $H_k\in [-\ve^{1/2},
\ve^{1/4}]$. We apply formula (4.1) with the change of the estimates
(4.2)\,--\,(4.6) as follows: when $x\in \cS_k([0,\ve^{1/4}])$ there
is no change in the estimates; when $x\in \cS_k([-\ve^{1/2},0])$ we
replace (4.2) and (4.3) by $K_1(x)=\dfrac{\pt}{\pt H_k}(K_1(x))=0$
and (4.4)\,--\,(4.6) remain the same. The function $K_4(x)=1$ and
$K_5(x)$ is then defined in a same way as in Lemma 4.2 with an
estimate $0<C_{34}\leq K_5(x)\leq C_{35}<\infty$ for $x\in
\cS_k([-\ve^{1/2},0])$. Again we let
$$u_k^{(3)}(H_k)=\int_{-\ve^{1/2}}^{H_k}\dfrac{K_6(\ve)-y}{K_4(y)}dy$$
with
$$K_6(\ve)=\left(\int_{-\ve^{1/2}}^{\ve^{1/4}}\dfrac{dy}{K_4(y)}\right)^{-1}
\left(\int_{-\ve^{1/2}}^{\ve^{1/4}}\dfrac{ydy}{K_4(y)}\right) \ .$$

It is then checked that $u_k(-\ve^{1/2})=u_k(\ve^{1/4})=0$ and
$K_6(\ve)\leq 2\ve^{1/4}$. The estimate (4.9) is still working for
$H_k\in [-\ve^{1/2},\ve^{1/4}]$. The estimates (4.10), (4.11) and
(4.12) are still working. Let $\zt([-\ve^{1/2},\ve^{1/4}])$ be the
first time when the process $X_t^\ve$ starting from a point $x\in
\gm$, first hits $\gmout$ or $\gmin$. A similar statement of (4.13)
is then obtained. We have

$$\lim\li_{\ve \da 0}\sup\li_{x\in \cS([-\ve^{1/2},\ve^{1/4}])}\E_x
\zt([-\ve^{1/2},\ve^{1/4}])=0 \ . \eqno(4.21)$$

The estimate of the hitting probability is a bit simpler. We
construct a barrier function $u_k^{(4)}$ similarly as in Lemma 4.3.
The function $K_7(H_k)$ is defined as in Lemma 4.3. But now we have
the property that $|K_7(H_k)|\leq C$ for $H_k\in [-\ve^{1/2},0]$. We
let
$$u_k^{(4)}(H_k)=1-\dfrac{\play{\int_{-\ve^{1/2}}^{H_k}\exp\left(-\int_{-\ve^{1/2}}^y
K_7(z)dz\right)dy}}{\play{\int_{-\ve^{1/2}}^{\ve^{1/4}}\exp\left(-\int_{-\ve^{1/2}}^y
K_7(z)dz\right)dy}} \ .$$ We have $u_k^{(4)}(-\ve^{1/2})=1$,
$u_k^{(4)}(\ve^{1/4})=0$. We have

$$C_{36}\ve^{1/2}\leq
\int_{-\ve^{1/2}}^0\exp\left(-\int_{-\ve^{1/2}}^y
K_7(z)dz\right)dy\leq C_{37}\ve^{1/2} \ , $$ and by the second
inequality in Lemma 4.5 we see that

$$C_{38}\ve^{1/2}\leq
\int_{0}^{\ve^{1/4}}\exp\left(-\int_{-\ve^{1/2}}^y K_7(z)dz\right)dy
\leq C_{39}\ve^{1/2} \ .$$

These estimates ensure that an analogue of (4.19) works, but the
lower bound is a positive constant, and hence situation is a bit
simpler. We have

$$1\geq \Prob_x(X_{\zt([-\ve^{1/2},\ve^{1/4}])}^\ve\in \gmin)\geq u_k^{(4)}(0)\geq C_{40}>0 \ . \eqno(4.22)$$
uniformly in $x\in\gm$. The results (4.21), (4.22) and Lemma 4.8,
combined with a similar analysis of (4.20) in Lemma 4.6, give the
statement of this Lemma. $\square$

\

\textbf{Lemma 4.10.} \textit{We have $$\lim\li_{\ve\da
0}\sup\li_{x\in \gmin} \E_x\sm=0 \ .$$}

\

\textbf{Proof.} This is a result in the same essence of Lemma 3.2
(formula (10)) of \cite{[DK2]}. $\square$

\

\textbf{Lemma 4.11.} \textit{The process $X_{\sm_n}^\ve$ satisfies
the Doeblin condition on $\gm$.}

\

\textbf{Proof.} For each fixed $\ve>0$ we have the ergodicity of the
process $X_t^\ve$. The Doeblin condition is satisfied for the
process $X_t^\ve$ in $[\cE]$ and each of these $[U_k]$'s for
$k=1,...,r$. As we have Lemmas 4.8, 4.9 and 4.10, we see that the
statement of this Lemma follows. $\square$

\

\textbf{Acknowledgement.} We would first like to thank the anonymous
referee for the very careful reading of the manuscript and for
pointing out some errors that we did not notice, as well as giving
various important suggestions. We would also like to thank
Konstantinos Spiliopoulos for reading the first version of this
paper and to thank Rongrong Wang for a pleasant discussion on
problems about matrix perturbations.

\

\end{document}